\documentclass[11pt,a4paper]{article}
\usepackage[T1]{fontenc}
\usepackage{lmodern}
\usepackage[margin=27mm]{geometry}
\usepackage{amsmath,amssymb,amsthm,mathtools}
\usepackage{microtype}
\usepackage{enumitem}
\usepackage{tikz}
\usepackage{mathrsfs}
\usepackage{xcolor}
\usepackage{algorithm}
\usepackage{algpseudocode}
\usepackage{float}
\usepackage{booktabs}
\usepackage[T1]{fontenc}
\usepackage{lmodern}
\usepackage[colorlinks=true,linkcolor=blue!45!black,citecolor=blue!45!black,urlcolor=blue!45!black]{hyperref}
\usetikzlibrary{positioning}
\numberwithin{equation}{section}
\newtheorem{theorem}{Theorem}[section]
\newtheorem{proposition}[theorem]{Proposition}
\newtheorem{lemma}[theorem]{Lemma}
\newtheorem{corollary}[theorem]{Corollary}
\theoremstyle{remark}
\newtheorem{remark}[theorem]{Remark}
\newcommand{\Pbb}{\mathbb P}
\newcommand{\Ebb}{\mathbb E}

\newcommand{\eps}{\varepsilon}

\hypersetup{pdftitle={2 neighborhood reconstruction at the polylogarithmic threshold}}
\title{The second reconstruction transition for random graphs from 2-neighborhoods}
\author{
\textsc{Yiming Chen}\thanks{
School of Mathematical Sciences, Peking University,
Beijing 100871, China.
Email: \texttt{ymchenmath@math.pku.edu.cn}
}
}
\date{}
\begin{document}
\maketitle

\begin{abstract}
We prove that, with high probability, the Erd\H{o}s--R\'enyi random graph
$G(n,p)$ is reconstructible up to isomorphism from the multiset of its
rooted induced 2-neighborhoods whenever
$p\geq n^{-3/4}(\log n)^3$.
Together with the non-reconstructibility result of Johnston, Kronenberg, Roberts and Scott~\cite{JohnstonKronenbergRobertsScott}, our result locates the second reconstruction transition at the scale $p=n^{-3/4}$ up to polylogarithmic factors, resolving a question they posed.
Our proof relies on a global rigidity argument exploiting 5-cycle constraints.
\end{abstract}
\par

\section{Introduction}
\label{sec:introduction}

The graph shotgun assembly problem asks whether a graph is uniquely determined, up to isomorphism, by the multiset of its rooted neighborhoods. Although each neighborhood has a distinguished root, the overlaps between different neighborhoods are not specified. The problem is motivated by DNA shotgun sequencing, where a long sequence is assembled from overlapping short reads and repeated segments can prevent unique reconstruction~\cite{MotahariBreslerTse,BreslerBreslerTse}. Mossel and Ross~\cite{MosselRoss} introduced graph shotgun assembly for random models and established general conditions for reconstruction. In this paper, we consider the reconstruction of Erd\H{o}s--R\'enyi random graphs $G(n,p)$ from their rooted radius 2-neighborhoods.

For a graph $G$ and a vertex $v\in V(G)$, let $B_2(G,v)=G[\{u\in V(G):d_G(u,v)\le 2\}]$ be subgraph induced by the vertices at distance at most $2$ from $v$. The $2$-neighborhood deck of $G$ is the multiset
\[
D_2(G)=\bigl\{\!\bigl\{[B_2(G,v),v]:v\in V(G)\bigr\}\!\bigr\},
\]
we say that $G$ is reconstructible if every graph $H$ with $D_2(H)=D_2(G)$ is isomorphic to $G$.

Gaudio and Mossel~\cite{GaudioMossel} studied reconstruction from $1$- and $2$-neighborhoods of Erd\H{o}s--R\'enyi random graphs. For $p=n^{-\alpha}$, they proved reconstructibility from $2$-neighborhoods when $0<\alpha<1/2$ or $1/2<\alpha<3/5$, and non-reconstructibility when $3/4<\alpha<1$. For $1$-neighborhoods, Huang and Tikhomirov~\cite{HuangTikhomirov} subsequently determined the second reconstruction transition at the scale $p=n^{-1/2}$, up to polylogarithmic factors. 

Johnston, Kronenberg, Roberts and Scott~\cite{JohnstonKronenbergRobertsScott} substantially improved the bounds for $2$-neighborhoods. They proved that, for some constant $\delta>0$, $G(n,p)$ is reconstructible with high probability when $p\geq n^{-2/3-\delta}$. They also showed that reconstruction fails with high probability when
\begin{equation}
p=\omega(n^{-5/4})
\qquad\text{and}\qquad
p\leq \frac13 n^{-3/4}(\log n)^{1/4}.
\label{eq}
\end{equation}

\begin{samepage}
Despite these advances, the polynomial scale of the second reconstruction transition for $2$-neighborhoods remained undetermined. Our main theorem closes this gap.

\begin{theorem}
\label{thm:main}
If \begin{equation}
 n^{-3/4}(\log n)^3\leq p\leq n^{-2/3-\gamma},
 \label{eq:new-density-range}
\end{equation}
where $\gamma>0$ is a sufficiently small fixed constant. Then $G(n,p)$ is reconstructible with high probability.
\end{theorem}
\end{samepage}

Together with~\eqref{eq}, Theorem~\ref{thm:main} locates the second reconstruction transition for $2$-neighborhoods at the scale $p=n^{-3/4}$, up to polylogarithmic factors.



\subsection{Related work}

\paragraph{Graph reconstruction.}

The classical reconstruction conjecture asks whether every graph with at least three vertices is determined by the multiset of subgraphs obtained by deleting a vertex~\cite{Kelly,Ulam}. M\"uller~\cite{Muller} proved that almost every graph is reconstructible, and Bollob\'as~\cite{BollobasReconstruction} proved that almost every graph has reconstruction number three. For further results, see ~\cite{Kelly,HararyPlantholt,AsilisChenHansenTeng}. 

%

A related line of work concerns canonical labeling in random graphs. Erd\H{o}s and R\'enyi~\cite{ErdosRenyiAsymmetric} showed that almost all graphs are asymmetric, while Babai, Erd\H{o}s and Selkow~\cite{BabaiErdosSelkow} and Babai and Ku\v{c}era~\cite{BabaiKucera} developed algorithms for canonical labeling. Bollob\'as~\cite{BollobasDistinguishing} used distance sequences to distinguish vertices, and Czajka and Pandurangan~\cite{CzajkaPandurangan} used the degrees of neighboring vertices to obtain canonical labels in sparse random graphs. Gaudio, R\'acz and Sridhar~\cite{GaudioRaczSridhar} showed that $3$-neighborhoods provide canonical labels near the connectivity threshold, and extended their analysis to randomly perturbed graphs. Verbitsky and Zhukovskii~\cite{VerbitskyZhukovskii} obtained canonical labels from counts of short walks. In the range~\eqref{eq:new-density-range}, the multisets of neighboring degrees distinguish all vertices with high probability~\cite[proof of Theorem~4.5]{MosselRoss}. 


\paragraph{Shotgun assembly.}

Mossel and Ross~\cite{MosselRoss} showed that the reconstruction radius for $G(n,\lambda/n)$ is of order $\log n$ for every fixed $\lambda>0$ with $\lambda\neq1$. Ding, Jiang and Ma~\cite{DingJiangMa} subsequently determined its sharp asymptotics for all fixed $\lambda>0$, including $\lambda=1$. More precisely, if $T,T'$ are independent Poisson--Galton--Watson trees with mean $\lambda$, the reconstruction threshold satisfies
\[
r_*\sim\frac{\log n}{\log\bigl((\lambda^2\mathbb P(T\cong_{\mathrm{root}}T'))^{-1}\bigr)}.
\]
They also gave a reconstruction algorithm that runs in polynomial time above the threshold. For random $d$-regular graphs with fixed $d\geq3$, Huang, Mossel, Sun, Zhang and Zhou~\cite{HuangMosselSunZhangZhou} determined the reconstruction radius up to an additive constant.

Beyond random graphs, Ding and Liu~\cite{DingLiuLattice} determined the sharp reconstruction threshold for independent random labels on lattice boxes and gave an efficient reconstruction algorithm. Adhikari and Chakraborty~\cite{AdhikariChakrabortyLM} studied shotgun assembly for random simplicial complexes in the Linial--Meshulam model. Related problems have also been studied for random jigsaw puzzles~\cite{BordenaveFeigeMossel,NenadovPfisterSteger,MartinssonJigsaw,BalisterBollobasNarayanan}, random colorings of the hypercube~\cite{PrzykuckiRobertsScott}, random pictures~\cite{NarayananYap}, and patterns on groups~\cite{RaymondBlandMcGoff}.

\paragraph{Random graph matching.}

Random graph matching seeks to recover the unknown vertex correspondence between two correlated random graphs. Ding, Ma, Wu and Xu~\cite{DingMaWuXuProfiles} developed efficient matching algorithms based on degree profiles, which can be computed from $2$-neighborhoods. Ding, Fei and Wang~\cite{DingFeiWangProfiles} extended this approach to random graphs with inhomogeneous edge probabilities. For sparse graphs, comparisons of neighborhood trees and tests for correlation between random trees have led to efficient algorithms for partial alignment~\cite{GanassaliMassoulieTree,GanassaliLelargeMassoulieCorrelation}. Mao, Wu, Xu and Yu~\cite{MaoWuXuYuChandeliers} developed an efficient matching algorithm based on counts of specially chosen rooted trees, called chandeliers. Wu, Xu and Yu~\cite{WuXuYuThresholds} studied thresholds for exact and partial recovery, while Ding and Du~\cite{DingDuRecovery} determined the sharp threshold for recovering a positive fraction of vertices in sparse random graphs. Ding and Li~\cite{DingLiGraph} obtained an efficient matching algorithm for any fixed positive edge correlation in a broad range of densities.  Related work also addresses detection, recovery, computational barriers, independent graph alignment, and reconstruction from graph projections~\cite{DingDuDetection,GanassaliMassoulieSemerjian,DuOptimalRecovery,DingDuLiLowDegree,ChenDingGongLiTransition,ChenDingGongLiDecoratedTrees,DuGongHuangAlignment,GongLiXuProjection}.

\subsection{Proof overview}

For $p$ in the range~\eqref{eq:new-density-range}, let $\mathcal C_{n,p}$ denote a class of graphs, to be specified later, such that
\[
\mathbb P(G(n,p)\in\mathcal C_{n,p})=1-o(1).
\]

Theorem~\ref{prop:approx} shows that, with high probability, every $H\in\mathcal C_{n,p}$ whose rooted $2$-neighborhoods are isomorphic to those of $G$ at every vertex satisfies
\[
|E(G)\triangle E(H)|
=o(n^2p).
\]
On the other hand, Theorem~\ref{thm:correction} provides a class $\mathcal T_{n,p}\subseteq\mathcal C_{n,p}$, with $\mathbb P(G\in\mathcal T_{n,p})=1-o(1)$. For some constant $c=c(\gamma)>0$, any two distinct graphs $G,H\in\mathcal T_{n,p}$ whose rooted $2$-neighborhoods are isomorphic at every vertex satisfy
\[
|E(G)\triangle E(H)|\geq cn^2p.
\]
Both bounds rely on the global rigidity imposed by the $5$-cycle
constraints. These constraints
restrict compatible edge pairings in the upper bound and rule out small
nonzero edge discrepancies in the lower bound. Combining the two bounds, we obtain a event $\mathcal U_{n,p}$. No two nonisomorphic graphs in $\mathcal U_{n,p}$ have the same deck. Finally, Lemma~\ref{lem:fibre} gives
\[
\mathbb P(G\text{ is not reconstructible})
\leq 2\mathbb P(G\notin\mathcal U_{n,p})=o(1),
\]
completing the proof.

\paragraph{Organization.}

Section~\ref{sec:prelim} introduces preliminaries. Section~\ref{sec:approx} establishes an upper bound on the symmetric difference of edge sets, while Section~\ref{sec:correction} provides a lower bound for nonzero discrepancies. Section~\ref{sec:fibres} completes the proof of the main theorem. Appendices~\ref{app:preliminaries} and~\ref{app:correction} contain the proofs of the preliminaries and Proposition~\ref{prop:moved}.

\section{Preliminaries}\label{sec:prelim}
Unless stated otherwise, we fix a
sufficiently small constant $\gamma>0$ and
assume~\eqref{eq:new-density-range}. Throughout the paper, $C,c>0$ denote constants
whose values may vary from line to line.
For a graph $F$, let $N_F(u)$ denote the neighborhood of $u$,
and let $P_j^F(x,y)$ denote the number of simple paths
of length $j$ between $x$ and $y$. We write $P_j^{F-u}(x,y)$ for the corresponding count
in the graph obtained from $F$ by deleting $u$. A 5-cycle in a graph $F$ is a subgraph isomorphic to $C_5$,
not necessarily induced. For $e\in E(F)$, let $c_5(F;e)$ denote the number
of 5-cycles containing $e$.

For each vertex $u\in V(F)$, define the 5-cycle link
$Q_u(F)$ to be the graph with vertex set $N_F(u)$,
where
\[
ab\in E(Q_u(F))
\quad\Longleftrightarrow\quad
P_3^{F-u}(a,b)>0.
\]
Every witnessing path $a-x-y-b$, together with $ua$ and $ub$, lies in
$B_2(F,u)$. 

\subsection{The degree fingerprints}

For a graph $F$ and $v\in V(F)$, define the degree fingerprints
of $v$ by
\[
\chi_F(v)
=\left(\deg_F(v),
\{\deg_F(w):w\in N_F(v)\}_{\mathrm{multiset}}\right).
\]
$\chi_F(v)$ is determined by the rooted
$2$-neighborhood $(B_2(F,v),v)$.

\begin{lemma}[Mossel--Ross~{\cite[proof of Theorem~4.5]{MosselRoss}}]
\label{lem:fingerprints}
Let $G\sim G(n,p)$, where
$np/(\log n)^2\to\infty$ and $p\le n^{-3/5}$.
Then, with high probability, 
$\chi_G(v)$, $v\in V(G)$, are pairwise distinct.
\end{lemma}

The use of degree profiles for canonical labeling
goes back to Czajka and Pandurangan~\cite{CzajkaPandurangan};
see also Gaudio, R\'acz and Sridhar~\cite{GaudioRaczSridhar}. Note that \eqref{eq:new-density-range}
satisfy the conditions of Lemma~\ref{lem:fingerprints}.
Consequently, with high probability, every graph $H$
with $D_2(H)=D_2(G)$ admits a unique bijection
$\phi:V(G)\to V(H)$ such that $\chi_G(v)=\chi_H(\phi(v))$ for all $v\in V(G).$
Lemma~\ref{lem:fingerprints} thus allows us to assign
a unique label to each vertex based on its degree profile. For graphs $F,H$ on $[n]$, write $F\equiv_2 H$ if
$B_2(F,u)\cong_{\mathrm{root}} B_2(H,u)$
for every $u\in[n]$.

\subsection{Uniform geometry of the random graph}

\eqref{eq:new-density-range} imply
$np\ge n^{1/4}(\log n)^3$, $n^2p^3\le n^{-3\gamma}$,
and $n^3p^4\ge(\log n)^{12}$.
In this section, we introduce some basic estimates for degrees,
short paths, and 5-cycles, together with bounds on the
local geometry of 5-cycle links.

\begin{lemma}\label{lem:basic}
For every fixed $M>0$, there exists a constant $C>0$
such that, with probability $1-O(n^{-M})$,
every degree lies in $[np/2,2np]$,
every codegree and every $P_3(x,y)$, $x\ne y$,
is at most $C$, and every vertex is in at most
$C$ triangles.
\end{lemma}

\begin{corollary}\label{lem:disjoint-fourpaths}
There exists a constant $C>0$ such that,
with high probability, for every pair of distinct
vertices $x,y$, every collection of $\kappa$ distinct
4-paths between $x$ and $y$ contains at least
$\lfloor \kappa/C\rfloor$ paths with pairwise disjoint
internal vertex sets.
\end{corollary}

\begin{lemma}\label{lem:fourpaths}
There exists a constant $C>0$ such that,
with high probability, for every pair of
distinct vertices $x,y$,
\begin{equation}
P_4(x,y)\le Cn^3p^4.
\label{eq:P4_bound}
\end{equation}
\end{lemma}

\begin{lemma}\label{lem:cycles}
There exist constants $c,C>0$ such that,
with high probability, every edge is contained
in between $cn^3p^4$ and $Cn^3p^4$ 5-cycles,
and every pair of distinct edges is contained
in at most $C$ common 5-cycles.
\end{lemma}

For $z\ne u$, define
\[
 \Gamma_u(z)=\{w\in N(u)\setminus\{z\}:P_3^{G-u}(z,w)>0\}.
\]
$\Gamma_u(v)\cap\Gamma_u(z)$ consists of neighbors of $u$ can be reached from both $v$ and $z$ via a 3-path that does not pass through $u$. The vertices $v,z$ themselves need not be adjacent to $u$.

\begin{figure}[H]
    \centering
    \includegraphics[width=0.72\textwidth]
    {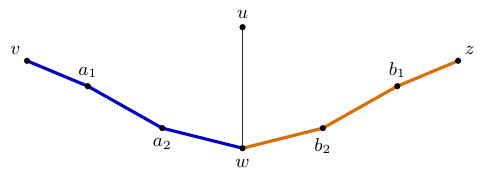}
    \caption{Two 3-paths meeting at $w$, witnessing
    $w\in\Gamma_u(v)\cap\Gamma_u(z)$.}
    \label{fig:path-overlap}
\end{figure}

\begin{lemma}\label{lem:overlap}
With high probability, for every triple of different vertices $u,v,z$,
%

\begin{equation}\label{eq:overlap}
 |\Gamma_u(v)\cap\Gamma_u(z)|
 \le C(n^5p^7+\log n).
\end{equation}
\end{lemma}

\begin{lemma}\label{lem:spectral}
Let $G\sim G(n,p)$, where $p>n^{-1}(\log n)^6$.
Let $A$ be the adjacency matrix of $G$.
With high probability, simultaneously for every $S\subseteq E(G)$,
\begin{equation}\label{eq:rope-count}
 \sum_{\substack{v_0,\ldots,v_4\in V(G)\\\text{all distinct}}}
 \mathbf1_{\{v_0v_1\in S\}}A_{v_1v_2}A_{v_2v_3}
 \mathbf1_{\{v_3v_4\in S\}}
 \le C\{np^2s^2+n^2p^2s\},
\end{equation}
where $s=|S|$ and $C>0$ is an absolute constant.
\end{lemma}


\begin{lemma}\label{lem:links}
For every fixed $\eta>0$ there is a fixed $\delta_\eta>0$ such that, with
high probability, simultaneously for every root $u$ and every $J\subseteq N(u)$,
\begin{equation}\label{eq:small-link-sets}
 |J|\le\delta_\eta np
 \quad\Longrightarrow\quad
 e(Q_u(G)[J])\le\eta n^3p^4|J|.
\end{equation}
\end{lemma}

For a graph $F$ and an edge set $B\subseteq E(F)$, let $J_F(B)$ be the number of pairs $(e,C)$ such that $e\in B$, $C$ is a simple 5-cycle of $F$ containing $e$, and $C$ contains another edge of $B$ with no endpoint in common with $e$. Each pair $(e,C)$ is counted once, regardless of how many such other edges occur in $C$.

\begin{proposition}\label{prop:single}
For every fixed $\varepsilon>0$, there is a fixed $\rho_\varepsilon>0$ such that, with high probability, simultaneously for every $B\subseteq E(G)$,
\begin{equation}\label{eq:small-edge-cycle-incidence}
 |B|\le\rho_\varepsilon n^2p
 \quad\Longrightarrow\quad
 J_G(B)\le\varepsilon n^3p^4|B|.
\end{equation}
\end{proposition}

\begin{remark}\label{rem:cycle-overlap}
Lemma~\ref{lem:links} and Proposition~\ref{prop:single}
consider different types of interactions among edges
in 5-cycles. For a fixed vertex $u$ and
$J\subseteq N(u)$, $e(Q_u(G)[J])$
counts pairs of edges $ua,ub$ with $a,b\in J$
that lie on a common 5-cycle, without counting
the multiplicity of such cycles. In contrast,
$J_G(B)$ counts pairs $(e,C)$ where $e\in B$
and $C$ is a 5-cycle containing $e$ and another edge of $B$ disjoint from $e$.
Thus, the former concerns edges sharing a vertex,
whereas the latter concerns two disjoint edges.
\end{remark}


\section{Proof of the upper bound of the symmetric difference}\label{sec:approx}

We first show that, with high probability, \(G\sim G(n,p)\) differs in only \(o(n^2p)\) edges from every typical graph \(H\) whose rooted \(2\)-neighborhoods are isomorphic to those of \(G\) at every vertex. We condition on the degree sequence and pass to the configuration model, representing the two graphs by perfect matchings of labelled half-edges. The common \(5\)-cycle links impose local restrictions: after fixing sparse partial matchings, a partner in one matching determines the vertex block of the corresponding partner in the other. We then condition on prescribed common pairs and count the possible alternating cycles formed by the remaining pairs. A union bound yields the \(o(n^2p)\) upper bound on the edge discrepancy.


Let $\mathcal C_{n,p}$ be the class of graphs satisfying the conclusions of
Lemmas~\ref{lem:basic}, \ref{lem:cycles} and~\ref{lem:overlap}. Then
$\Pbb(G(n,p)\in\mathcal C_{n,p})=1-o(1)$.

\begin{theorem}\label{prop:approx}
With probability $1-o(1)$, every $H\in\mathcal C_{n,p}$ on $[n]$ such that
$G\equiv_2 H$ satisfies
\begin{equation}
|E(G)\mathbin\triangle E(H)|
<\frac{2|E(G)|}{(\log n)^{1/4}}.
\label{eq:symdiff-bound}
\end{equation}
\end{theorem}

\begin{proof}


Fix a degree sequence $(r_u)_{u\in[n]}$
satisfying $np/2\le r_u\le 2np$.
For each vertex $u$, introduce a block of labelled half-edges $B_u=\{(u,1),\ldots,(u,r_u)\}.$ Let $N=\sum_{u\in[n]}r_u$ be the total number of
half-edges, and $m=N/2$ be the number of edges.
Choose a uniformly random perfect matching $P$
of these $N$ half-edges.
Each matched pair $\{s,t\}$ gives an edge between
the vertices to which $s$ and $t$ belong.
Write $F=\operatorname{proj}(P)$ for the resulting
multigraph on $[n]$, which may contain loops and
multiple edges.

Every simple graph with degree sequence $(r_u)$
has exactly $\prod_{u\in[n]}r_u!$ realizations
as a pairing. Since $G(n,p)$ conditioned on its
degree sequence is uniform over all simple graphs
with that sequence, we have
\begin{equation}\label{eq:sec3-configuration-law}
\mathcal L\bigl(F
\mid F\text{ is simple}\bigr)
=
\mathcal L\bigl(G(n,p)
\mid \deg(G)=(r_u)_{u\in[n]}\bigr).
\end{equation}
Thus, estimates for random pairings can be transferred
to $G(n,p)$ by conditioning on simplicity.

For each half-edge $s\in B_u$, write
$\operatorname{root}(s)=u$ and let $P(s)$ denote its mate. When $F$ is simple, the map $\nu_P(s):=\operatorname{root}(P(s)),$ $s\in B_u,$ is a bijection. Define the matrix $D_P$, indexed by the
half-edges, by $D_P(s,t)=
\mathbf 1_{\{\nu_P(s)\nu_P(t)\in E(Q_u(F))\}},$ $s,t\in B_u,$ and set $D_P(s,t)=0$ for half-edges in
different blocks. Thus, $D_P$ is the block diagonal adjacency matrix
of the $Q_u(F)$, with vertices indexed
by the corresponding half-edges.

For any pairing $P$ with simple projection $F$
and any graph $H$ satisfying $F\equiv_2 H$,
choose rooted isomorphisms
$f_u:B_2(H,u)\to B_2(F,u)$ for all $u\in[n]$. This construction is closely related to the
oriented edge correspondence in
\cite[Section~6]{HuangTikhomirov}.
For each edge $ux\in E(H)$, assign its end at $u$
to the unique half-edge $s\in B_u$ satisfying
$\nu_P(s)=f_u(x)$.
Pairing the half-edges assigned to the two ends
of each edge of $H$ yields a perfect matching
$Q$ with $\operatorname{proj}(Q)=H$.
Since the rooted isomorphisms preserve the
$5$-cycle links, we have
\begin{equation}\label{eq:sec3-common-link}
D_P=D_Q.
\end{equation}
We call $(P,Q)$ compatible if both projections
belong to $\mathcal C_{n,p}$ and
\eqref{eq:sec3-common-link} holds. 

We now extract small partial matchings
that determine the vertex block of a mate in one pairing
once its mate in the other pairing is specified. More precisely, for every compatible pair $(P,Q)$,
we claim that there exist partial matchings
$P_0\subseteq P$ and $Q_0\subseteq Q$ such that
\begin{equation}\label{eq:sec3-seed-size}
|P_0|+|Q_0|\le\frac{3m}{\sqrt{\log n}},
\end{equation}
and, for every half-edge $s$, the data
$P_0,Q_0,Q(s)$ determine $\nu_P(s)$, while
$P_0,Q_0,P(s)$ determine $\nu_Q(s)$.

To prove this, we fix $(P,Q)$, and
Independently retain each pair of $P$ and each
pair of $Q$ with probability $1/\sqrt{\log n}$,
obtaining $P_0$ and $Q_0$. Fix a half-edge $s\in B_u$.
By Lemma~\ref{lem:cycles} and $D_P=D_Q$,
there are at least $cn^3p^4/C$ half-edges
$t\in B_u$ such that
$$D_P(s,t)=D_Q(s,t)=1.$$

For each such $t$, choose a witnessing $5$-cycle
in each projection. Deleting the edge corresponding
to $s$ from each cycle leaves four pairs in $P$
and four in $Q$. We call a completion successful if all eight
pairs are retained in $P_0,Q_0$.

\begin{figure}[H]
    \centering
    \includegraphics[width=0.72\textwidth]
    {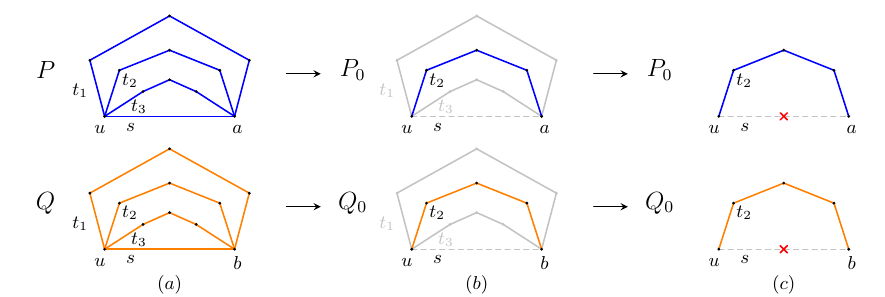}
    \caption{Construction of sparse partial matchings.
(a) Each projection contains $\Theta(n^3p^4)$
$5$-cycles through the edge at $s$.
(b) Independent thinning yields $P_0,Q_0$,
retaining a witnessing completion in each.
(c) Deleting the edge at $s$ leaves four pairs
in each of $P_0,Q_0$.}
    \label{fig:2}
\end{figure}

Each pair in $P$ or $Q$ appears in at most $C$
chosen completions, since distinct half-edges $t$
yield distinct $5$-cycles containing that pair and the edge at $s$.
Thus, the intersection graph of these completions has
maximum degree at most $8(C-1)$.
A greedy selection yields $\ell_s\ge cn^3p^4$
completions with pairwise disjoint pairs.

Let $Z_s$ count the eight pair completions retained in $P_0,Q_0$. By independence, conditional on $P,Q$, $Z_s\sim\operatorname{Bin}
\bigl(\ell_s,(\log n)^{-4}\bigr),$ and
\[
\Pbb\left(
Z_s<\frac{\ell_s}{2(\log n)^4}
\,\middle|\,P,Q
\right)
\le
\exp\left\{-\frac{\ell_s}{8(\log n)^4}\right\}.
\]

Since $n^3p^4/(\log n)^4\ge(\log n)^8$ and
$N=O(n^2p)$, a union bound yields
\begin{equation}
\Pbb\left(
\min_{s\in\bigcup_{u\in[n]}B_u} Z_s
\ge \frac{cn^3p^4}{2(\log n)^4}
\,\middle|\,P,Q
\right)=1-o(1).
\label{eq:min-Z-lower-bound}
\end{equation}

Moreover, $|P_0|+|Q_0|
\sim\operatorname{Bin}
\left(2m,\frac1{\sqrt{\log n}}\right),$ and
\[
\Pbb\left(
|P_0|+|Q_0|>
\frac{3m}{\sqrt{\log n}}
\,\middle|\,P,Q
\right)
\le
\exp\left\{-\frac{m}{5\sqrt{\log n}}\right\}
=o(1).
\]

Hence, every compatible pair $(P,Q)$ admits partial matchings
$P_0,Q_0$ satisfying \eqref{eq:sec3-seed-size},
and \eqref{eq:min-Z-lower-bound} holds
for all half-edges.

Now fix $P_0,Q_0$ as above and $s\in B_u$.
We show that $P_0,Q_0,Q(s)$ determine $\nu_P(s)$.
The key is to show that the true endpoint has
far more witnesses than any other candidate.
Write $b=\nu_Q(s)$ and let $R_s=\left\{r\in B_u:
\{r,Q(r)\}\in Q_0,\quad
P_3^{\operatorname{proj}(Q_0)-u}(\nu_Q(r),b)>0
\right\}.$ For each $r\in R_s$, the edges at $s$ and $r$
lie on a common $5$-cycle in $H$. Hence
$D_P(s,r)=D_Q(s,r)=1$ for every $r\in R_s$.
Writing $a=\nu_P(s)$, we obtain
\begin{equation}\label{eq:sec3-port-inclusion}
\nu_P(R_s)\subseteq\Gamma_u^F(a).
\end{equation}

For $w\ne u$, define $S_s(w)=\left|\left\{r\in R_s:
\{r,P(r)\}\in P_0,\quad
P_3^{\operatorname{proj}(P_0)-u}(\nu_P(r),w)>0
\right\}\right|.$ By the \eqref{eq:min-Z-lower-bound},
\[
S_s(a)\ge Z_s\ge
\frac{cn^3p^4}{2(\log n)^4}.
\]

On the other hand, each $r$ counted by $S_s(w)$
satisfies $\nu_P(r)\in\Gamma_u^F(w)$.
By \eqref{eq:sec3-port-inclusion} and the injectivity
of $\nu_P$ on $B_u$, for $w\notin\{u,a\}$,
\[
\begin{aligned}
S_s(w)
&\le |\Gamma_u^F(a)\cap\Gamma_u^F(w)|\\
&\le C(n^5p^7+\log n)
=o\left(\frac{n^3p^4}{(\log n)^4}\right),
\end{aligned}
\]
where we used $n^2p^3\le n^{-3\gamma}$ and
$n^3p^4\ge(\log n)^{12}$.

\begin{figure}[H]
    \centering
    \includegraphics[width=0.60\textwidth]{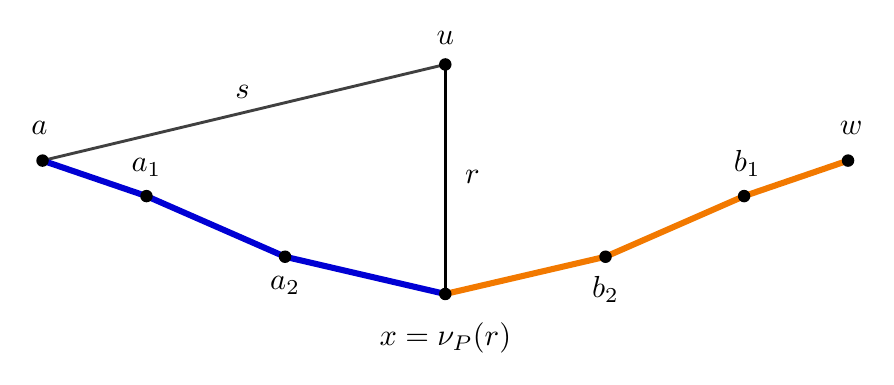}
    \caption{The edge at $s$, with true endpoint $a$, lies in
at least $cn^3p^4$ $5$-cycles, whereas each wrong
candidate $w\notin\{u,a\}$ is constrained by
$x=\nu_P(r)\in\Gamma_u^F(a)\cap\Gamma_u^F(w)$.}
    \label{fig:3}
\end{figure}

Thus, for all sufficiently large $n$,
\begin{equation}\label{eq:sec3-endpoint}
\nu_P(s)=
\underset{w\in[n]\setminus\{u\}}{\operatorname{argmax}}
S_s(w).
\end{equation}
By symmetry, fixing a mate in either pairing determines
the vertex block of its mate in the other.

We now bound the probability that $P$ admits a
compatible pairing $Q$ with many differing pairs
by counting alternating cycles.
Consider $k=|P\setminus Q|$.
Since both pairings contain $m$ pairs,
$|Q\setminus P|=k$ and $|P\cap Q|=m-k$.
Removing the common pairs leaves $2k$ half-edges
forming $h$ disjoint alternating cycles
(see \cite[Proof of Theorem~10.7]{KorteVygenCO}).
Each cycle contains at least two pairs from
each pairing, so $h\le k/2$.

Fix $k,h$, a collection of $m-k$ disjoint pairs
as a candidate for $P\cap Q$, and partial matchings
$P_0,Q_0$ satisfying \eqref{eq:sec3-seed-size}.
Also fix integers $\ell_1,\ldots,\ell_h\ge 2$
with $\sum_{i=1}^h\ell_i=k$, specifying the
numbers of $P$-pairs in the alternating cycles,
ordered by their least half-edges.

Let $\mathcal A$ be the event that there exists
a compatible $Q$ such that $P\cap Q$
equals the fixed collection of $m-k$ pairs,
$P_0\subseteq P$, $Q_0\subseteq Q$,
and the $i$th alternating cycle contains
$\ell_i$ pairs of $P$. We further require that both \eqref{eq:sec3-endpoint} and its counterpart with \(P\) and \(Q\) interchanged hold for every half-edge \(s\).

To count the admissible completions, traverse
the cycles in the prescribed order, starting
at the least unused half-edge and following
its $P$-pair first.
For a cycle containing $\ell$ pairs of $P$,
the first mate has at most $N$ choices,
each of the next $2\ell-2$ nonclosing steps
has at most $2np$ choices by the endpoint
restrictions, and the final pair is forced.
Thus, there are at most
$N^h(2np)^{2k-2h}$ admissible completions.

The probability that $P$ contains the fixed
common pairs is $\frac{(2k-1)!!}{(N-1)!!},$ and the remaining $2k$
half-edges are uniformly paired. Consequently,
\[
\begin{aligned}
\Pbb(\mathcal A)
&\le
\frac{(2k-1)!!}{(N-1)!!}
\cdot
\frac{N^h(2np)^{2k-2h}}{(2k-1)!!}
=
\frac{N^h(2np)^{2k-2h}}{(N-1)!!}.
\end{aligned}
\]

We now take a union bound over all choices.
The half-edges involved in $P\cap Q$, $P_0$, and $Q_0$
can be chosen in at most $2^{3N}$ ways.
Their pairings can then be specified by successively
matching the least unpaired half-edge with at most
$N$ choices at each step, giving at most
\[
N^{m-k+|P_0|+|Q_0|}
\le N^{m-k+3m/\sqrt{\log n}}
\]
possibilities, and the numbers of $P$-pairs in the $h$ cycles
can be chosen in at most
$\binom{k-1}{h-1}\le 2^k$ ways.
Consequently, using $(N-1)!!=\exp\{m\log N-m+O(\log N)\}$,
$N\ge n^2p/2$, $h\le k/2$, and
$p\le n^{-2/3-\gamma}$, we obtain
\begin{equation}\label{eq:sec3-pairing-exception}
\begin{aligned}
&\Pbb\left(
\exists Q:\ (P,Q)\text{ compatible},\
|P\setminus Q|\ge\frac{m}{(\log n)^{1/4}}
\right)\\
&\quad\le
\frac{2^{3N}N^{3m/\sqrt{\log n}}}{(N-1)!!}
\sum_{k=\lceil m/(\log n)^{1/4}\rceil}^{m}
\sum_{h=1}^{\lfloor k/2\rfloor}
2^kN^{m-k+h}(2np)^{2k-2h}\\
&\quad\le
\exp\{O(m\sqrt{\log n})\}
\sum_{k=\lceil m/(\log n)^{1/4}\rceil}^{m}
\sum_{h=1}^{\lfloor k/2\rfloor}
\left(\frac{4n^2p^2}{N}\right)^{k-h}\\
&\quad\le
m^2\exp\left\{
O(m\sqrt{\log n})
-\frac{m}{2(\log n)^{1/4}}
\log\frac1{8p}
\right\}\\
&\quad\le
\exp\{-cm(\log n)^{3/4}\},
\end{aligned}
\end{equation}
for some fixed $c>0$, uniformly over all
admissible degree sequences.

%
%

We now transfer the \eqref{eq:sec3-pairing-exception} to simple graphs.
Since $\frac{(\max_u r_u)^3}{N}=O(np^2)=o(1),$ Theorem~5.2 of~\cite{McKayWormald} gives, uniformly
over the prescribed degree sequences,
\begin{equation}\label{eq:sec3-simplicity}
\Pbb(F\text{ is simple})
\ge \exp\{-O(n^2p^2)\}.
\end{equation}
Combining \eqref{eq:sec3-simplicity} with \eqref{eq:sec3-pairing-exception}
and $n^2p^2=o(m(\log n)^{3/4})$, we obtain
\begin{equation}\label{eq:sec3-conditional-simplicity}
\begin{aligned}
&\Pbb\left(
\exists Q:\ (P,Q)\text{ compatible},\
|P\setminus Q|\ge\frac{m}{(\log n)^{1/4}}
\,\middle|\,F\text{ is simple}
\right)\\
&\qquad\le
\exp\{-cm(\log n)^{3/4}+O(n^2p^2)\}
=o(1).
\end{aligned}
\end{equation}

For simple projections $F,H$, distinct pairs
in $P\cap Q$ yield distinct common edges. Hence
\[
\begin{aligned}
|P\setminus Q|
&=m-|P\cap Q|\\
&\ge m-|E(F)\cap E(H)|\\
&=\frac12|E(F)\triangle E(H)|.
\end{aligned}
\]
Thus, if $F\in\mathcal C_{n,p}$ admits a competitor
$H\in\mathcal C_{n,p}$ with $F\equiv_2H$ and
\[
|E(F)\triangle E(H)|
\ge\frac{2m}{(\log n)^{1/4}},
\]
then every $P$ projecting to $F$
admits a compatible $Q$ satisfying $|P\setminus Q|\ge\frac{m}{(\log n)^{1/4}}.$ 
By \eqref{eq:sec3-configuration-law} and
\eqref{eq:sec3-conditional-simplicity},
averaging over the degree sequences yields
\[
\begin{aligned}
&\Pbb\left(
\exists H\in\mathcal C_{n,p}:
G\equiv_2H,\
|E(G)\triangle E(H)|
\ge\frac{2|E(G)|}{(\log n)^{1/4}}
\right)\\
&\qquad\le
\Pbb(G\notin\mathcal C_{n,p})+o(1)
=o(1).
\end{aligned}
\]

\end{proof}

For completeness, we provide an exponential time algorithm for
reconstruction from a prescribed deck of rooted $2$-neighborhoods.
For a graph $H$ on $[n]$, let $\operatorname{Can}(H)$ denote the
labeled graph whose adjacency matrix is lexicographically smallest
among all relabelings of $H$.

Algorithm~\ref{alg:exact-reconstruction} enumerates half-edge pairings
and checks whether the resulting simple graphs realize the input deck.
The $5$-cycle-link comparison is an optional pruning step, the final
rooted-isomorphism test guarantees correctness.

\begin{algorithm}[H]
\caption{Reconstruction algorithm from a rooted $2$-neighborhood deck}
\label{alg:exact-reconstruction}
\begin{algorithmic}[1]
\Require $C_i=[K_i,o_i]$, $i\in[n]$
\Ensure $\operatorname{Can}(H)$ or a status flag
\State $r_i\gets\deg_{K_i}(o_i)$ for $i\in[n]$
\If{$\sum_{i=1}^n r_i$ is odd}
    \State \Return \textsf{InvalidDeck}
\EndIf
\State $B_i\gets\{(i,j):1\le j\le r_i\}$ for $i\in[n]$
\State $\mathscr R\gets\varnothing$
\ForAll{perfect matchings $P$ of $\bigcup_{i=1}^n B_i$}
    \State $H\gets\operatorname{proj}(P)$
    \If{$H$ is not simple}
        \State \textbf{continue}
    \EndIf
    \If{$Q_i(H)\not\cong Q_{o_i}(K_i)$ for some $i\in[n]$}
        \State \textbf{continue}
    \EndIf
    \If{$(B_2(H,i),i)\not\cong_{\mathrm{root}}(K_i,o_i)$
         for some $i\in[n]$}
        \State \textbf{continue}
    \EndIf
    \State $\mathscr R\gets\mathscr R\cup\{\operatorname{Can}(H)\}$
    \If{$|\mathscr R|=2$}
        \State \Return \textsf{Ambiguous}
    \EndIf
\EndFor
\If{$\mathscr R=\varnothing$}
    \State \Return \textsf{NoCompatibleGraph}
\EndIf
\State \Return the unique element of $\mathscr R$
\end{algorithmic}
\end{algorithm}

\section{Lower bound for small discrepancy}\label{sec:correction}

In this section, we establish the
complementary lower bound to Theorem~\ref{prop:approx}: any two distinct typical graphs whose
rooted 2-neighborhoods are isomorphic at every vertex must differ
in $\Omega(n^2p)$ edges. The proof is by contradiction. Suppose that two distinct
graphs $F,H$ satisfy $F\equiv_2 H$ with a sufficiently
small edge discrepancy. Proposition~\ref{prop:moved}
gives a lower bound on the loss of 5-cycle links
caused by moving neighbour labels. Proposition~\ref{prop:single} then gives a strictly smaller upper
bound on 5-cycle incidences that contradicts
the lower bound.

For $T\subseteq V(F)\setminus\{u\}$, let $\mathcal L_u^F(T)$ be the graph on $T$ joining distinct $x,y$ whenever $P_3^{F-u}(x,y)>0$. $\mathcal L_u^F(T)$ generalizes $Q_u(F)$, since
$Q_u(F)=\mathcal L_u^F(N_F(u))$.

\begin{proposition}\label{prop:moved}
There exists a constant $c>0$ such that, with high
probability, simultaneously for every root $u$,
every $T\subseteq[n]\setminus\{u\}$, and every bijection
$\phi:T\to N_G(u)$,
\begin{equation}\label{eq:moved-loss}
\left|
\phi^{-1}\big(E(Q_u(G))\big)
\setminus E(\mathcal L_u^G(T))
\right|
\ge cn^3p^4
|\{v\in T:\phi(v)\ne v\}|.
\end{equation}

\end{proposition}

Let $\mathcal T_{n,p}\subseteq\mathcal C_{n,p}$ be the
class of graphs satisfying \eqref{eq:moved-loss}, \eqref{eq:small-link-sets}, and \eqref{eq:small-edge-cycle-incidence}.
These properties are invariant under relabelling
and hold with high probability. Consequently,
\[
\mathbb P(G(n,p)\in\mathcal T_{n,p})=1-o(1).
\]

\begin{theorem}\label{thm:correction}
There exists a constant $c>0$, depending only on
$\gamma$, such that any two distinct graphs
$F,H\in\mathcal T_{n,p}$ satisfying $F\equiv_2 H$
obey
\[
|E(F)\triangle E(H)|\ge cn^2p
\]
for all sufficiently large $n$.
\end{theorem}

\begin{proof}

Fix distinct graphs $F,H\in\mathcal T_{n,p}$
satisfying $F\equiv_2 H$.
All constants in the proof are fixed independently
of $F$ and $H$. For each $u\in[n]$, fix a rooted isomorphism $f_u:B_2(H,u)\to B_2(F,u)$ as before.
Write $B=E(H)\setminus E(F)$ and $k'=|B|$.
Since $F\ne H$ and the degrees agree at every root,
we have $|E(F)\setminus E(H)|=k'>0$. Let $D_u=\{v\in N_H(u):f_u(v)\ne v\}.$ Since $f_u$ maps $N_H(u)$ bijectively onto $N_F(u)$,
we have $N_H(u)\setminus N_F(u)\subseteq D_u,$ and $\deg_B(u)\le |D_u|.$


We first bound the total number of moved neighbour labels. By Proposition~\ref{prop:moved}, there are at least
$c|D_u|n^3p^4$ unordered pairs
$\{x,y\}\subseteq N_H(u)$ satisfying $f_u(x)f_u(y)\in E(Q_u(F)),$ and $P_3^{F-u}(x,y)=0.$ Since $f_u$ is a rooted isomorphism, each such pair
is joined by a 3-path in $H-u$, which must contain
an edge $e\in B$. Together with the edges $ux$ and $uy$,
this path forms a 5-cycle $C$ in $H$, with $e$
not incident with $u$.

Distinct pairs at the same root yield distinct cycles,
and each pair $(e,C)$ can be associated with at most
three roots, namely the vertices of $C$ outside $e$.
Summing over all roots and applying
Lemma~\ref{lem:cycles}, we obtain
\begin{equation}\label{eq:total-movement}
\sum_{u\in[n]}|D_u|\le Ck'.
\end{equation}

Choose a sufficiently small fixed constant $\eta>0$,
and let $\delta>0$ be as
in~\eqref{eq:small-link-sets}.
For every $u\in[n]$ and $D\subseteq N_H(u)$
with $|D|\le\delta np$, \eqref{eq:small-link-sets}
gives
\[
e(Q_u(H)[D])\le\eta n^3p^4|D|.
\]

Now let $U=\{u\in[n]:|D_u|>\delta np\}$, $E_H(U)$ be the set of $H$-edges incident to $U$, and define $B^+=B\cup E_H(U)$.
By~\eqref{eq:total-movement} and $\max_{u\in[n]}\deg_H(u)\le 2np.$ Thus,
\begin{equation}\label{eq:correction-augmented-set}
|B^+|\le k'+2np|U|\le Ck'.
\end{equation}

For each $u\in U$ and each pair $(e,C)$ arising
from the construction at $u$ as above,
at least one of the two edges of $C$ incident
with $u$ is disjoint from $e$.
Since this edge belongs to $E_H(U)$ and $e\in B$,
the pair $(e,C)$ is counted by $J_H(B^+)$.
As each pair $(e,C)$ can arise from at most three roots,
we obtain
\begin{equation}\label{eq:correction-root-movement}
c n^3p^4\sum_{u\in U}|D_u|\le 3J_H(B^+).
\end{equation}

\begin{figure}[H]
    \centering
    \includegraphics[width=0.4\linewidth]
    {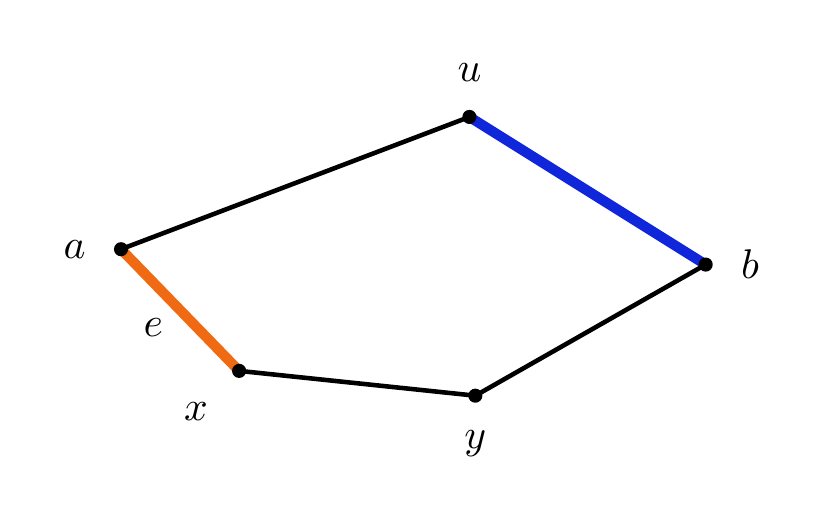}
    \caption{A 5-cycle $C$ in $H$ through $u\in U$.
    The edge $e\in B$ (orange) is disjoint from
    the edge $ub\in E_H(U)$ (blue), so $(e,C)$
    is counted by $J_H(B^+)$.}
    \label{fig:correction-edges}
\end{figure}

We now count all pairs $(e,C)$ with $e\in B$
and $C$ a 5-cycle of $H$ containing $e$.
By Lemma~\ref{lem:cycles}, there are at least
$ck'n^3p^4$ such pairs.
Those for which $e$ is incident to $U$ number at most
\[
Cn^3p^4\sum_{u\in U}\deg_B(u)
\le Cn^3p^4\sum_{u\in U}|D_u|
\le CJ_H(B^+),
\]
where the last inequality follows from
\eqref{eq:correction-root-movement}.
Among the remaining pairs, those for which $C$
contains another edge of $B$ disjoint from $e$
number at most $J_H(B)\le J_H(B^+)$.

Consider a remaining pair $(e,C)$. Both endpoints of $e$
lie outside $U$, and every other edge of $B$ on $C$
shares a vertex with $e$. If such an edge exists, choose its common endpoint $u$
with $e$. Then both edges of $C$ incident with $u$
belong to $B$, so both neighbours of $u$ on $C$
lie in $D_u$. Otherwise, orient $e=uy$ by a fixed rule, and let $w$
be the other neighbour of $u$ on $C$.
Since $y\in D_u$, the same conclusion holds if
$w\in D_u$.

In either case, assign $(e,C)$ to one such pair $(u,C)$
by a fixed rule. Each $(u,C)$ receives at most two
assignments, since $e$ is incident with $u$.
By Lemma~\ref{lem:basic}, each pair of neighbours
of $u$ lies on at most $C$ 5-cycles through $u$.
Hence, using~\eqref{eq:small-link-sets}
and~\eqref{eq:total-movement}, the number of
pairs assigned in this way is at most
\[
C\sum_{u\notin U}e(Q_u(H)[D_u])
\le C\eta n^3p^4\sum_{u\notin U}|D_u|
\le C\eta k'n^3p^4.
\]

\begin{figure}[H]
    \centering
    \includegraphics[width=0.80\linewidth]
    {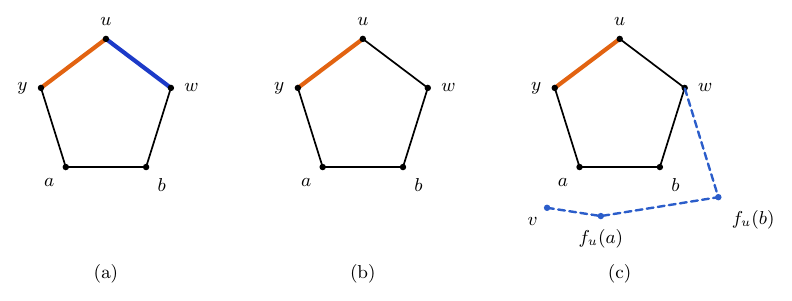}
    \caption{Three configurations for a remaining pair $(e,C)$
with $u\notin U$ and $e=uy\in B$ (orange).
(a) The cycle $C$ contains another $B$-edge $uw$ (blue).
(b) The edge $e$ is the only $B$-edge of $C$ and $w\in D_u$.
(c) The edge $e$ is the only $B$-edge of $C$ and $w\notin D_u$.
The dashed path in (c) represents the image under $f_u$
of the common-edge path $y-a-b-w$, joining $v=f_u(y)$
to $w$ in $F-u$.}
\label{fig:correction-cases}
    \label{fig:correction-cases}
\end{figure}

It remains to consider pairs $(e,C)$ for which
$e=uy$ is the only edge of $B$ on $C$ and
$w\notin D_u$. Then $w\in N_F(u)$ and $f_u(w)=w$.
Let $v=f_u(y)\ne y$, since $uy\notin E(F)$.
The 3-path from $y$ to $w$ along $C-u$ consists
of common edges and hence lies in $F-u$.
Since $C\subseteq B_2(H,u)$, applying $f_u$
to this path also gives a 3-path from $v$ to $w$
in $F-u$. Thus $w\in\Gamma_u^F(y)\cap\Gamma_u^F(v).$ By Lemma~\ref{lem:overlap}, there are uniformly
$o(n^3p^4)$ possible choices for $w$.
For each such $w$, Lemma~\ref{lem:basic} gives at most
$C$ 5-cycles in $H$ containing both $uy$ and $uw$.

Combining all the cases, we obtain
\begin{equation}\label{eq:correction-incidence-count}
ck'n^3p^4
\le CJ_H(B^+)+C\eta k'n^3p^4+o(k'n^3p^4).
\end{equation}

By our choice of $\eta$, the second term
on the right-hand side of
\eqref{eq:correction-incidence-count}
is at most $(c/4)k'n^3p^4$.
Next, fix $\varepsilon>0$ sufficiently small that
$C^2\varepsilon\le c/4$, and let $\rho_\varepsilon>0$
be given by Proposition~\ref{prop:single}.
Suppose that $0<k'\le \rho_\varepsilon n^2p/C$.
By~\eqref{eq:correction-augmented-set},
we have $|B^+|\le \rho_\varepsilon n^2p$.
Hence Proposition~\ref{prop:single} gives
\[
J_H(B^+)
\le\varepsilon n^3p^4|B^+|
\le C\varepsilon k'n^3p^4.
\]
Substituting into~\eqref{eq:correction-incidence-count}
yields
\[
ck'n^3p^4
\le \left(\frac{c}{2}+o(1)\right)k'n^3p^4,
\]
a contradiction for all sufficiently large $n$.
Therefore $k'>\rho_\varepsilon n^2p/C$, and
\[
|E(F)\triangle E(H)|=2k'=\Omega(n^2p).
\]

\end{proof}

\section{From typical graphs to arbitrary competitors}\label{sec:fibres}
We complete the proof by using the distribution of
$G(n,p)$ among graphs with the same deck to extend
uniqueness within a high-probability class to
reconstructibility with high probability.

Let $\mathcal U_{n,p}$ consist of all graphs
$F\in\mathcal T_{n,p}$ with pairwise distinct
degree fingerprints such that every
$H\in\mathcal C_{n,p}$ satisfying $F\equiv_2 H$
obeys \eqref{eq:symdiff-bound}. By Lemma~\ref{lem:fingerprints} and Theorem~\ref{prop:approx}, 
the class $\mathcal U_{n,p}$ is invariant under
relabelling and satisfies $\Pbb(G(n,p)\in\mathcal U_{n,p})=1-o(1).$

\begin{lemma}\label{lem:fibre}
Suppose no two nonisomorphic graphs in
$\mathcal U_{n,p}$ have the same
2-neighborhood deck. Then
\begin{equation}\label{eq:fibre-bound}
\Pbb\bigl(G(n,p)\text{ is not reconstructible}\bigr)
\le 2\Pbb\bigl(G(n,p)\notin\mathcal U_{n,p}\bigr).
\end{equation}
\end{lemma}

\begin{proof}
Fix a deck $\mathcal D$ represented in $\mathcal U_{n,p}$.
Since the deck determines the degree fingerprints,
every graph with deck $\mathcal D$ has a trivial
automorphism group. Hence each isomorphism class
has exactly $n!$ labelled representatives.

Moreover, let $\deg_{\mathrm{root}}(\mathsf d)$ denote
the degree of the root in a card $\mathsf d$.
Then $m=\frac12\sum_{\mathsf d\in\mathcal D}
\deg_{\mathrm{root}}(\mathsf d),$
\begin{equation}
\Pbb(G(n,p)\cong F)
=n!p^m(1-p)^{\binom n2-m}
\label{eq:class-prob}
\end{equation}
for every graph $F$ with deck $\mathcal D$.

If more than one isomorphism class has deck $\mathcal D$,
exactly one lies in $\mathcal U_{n,p}$,
while all the others lie outside it. By~\eqref{eq:class-prob}, all isomorphism classes
with deck $\mathcal D$ have equal probability.
Summing over all such decks, we obtain
\[
\begin{aligned}
&\Pbb(G\in\mathcal U_{n,p},
G\text{ is not reconstructible})\\
&=\sum_{\mathcal D}
\Pbb(G\in\mathcal U_{n,p},
G\text{ is not reconstructible},
\mathcal D_2(G)=\mathcal D)\\
&\le\sum_{\mathcal D}
\Pbb(G\notin\mathcal U_{n,p},
\mathcal D_2(G)=\mathcal D)\\
&=\Pbb(G\notin\mathcal U_{n,p}).
\end{aligned}
\]
Finally, 
\[
\begin{aligned}
\Pbb(G\text{ is not reconstructible})
&\le \Pbb(G\in\mathcal U_{n,p},
G\text{ is not reconstructible})
+\Pbb(G\notin\mathcal U_{n,p})\\
&\le 2\Pbb(G\notin\mathcal U_{n,p}).
\end{aligned}
\]
\end{proof}

\begin{proof}[Proof of Theorem~\ref{thm:main}]
Fix a sufficiently small $0<\gamma<1/12$,
and first consider $p$ in the range
\eqref{eq:new-density-range}.

Suppose two nonisomorphic graphs
$F,H\in\mathcal U_{n,p}$ have the same deck.
Matching vertices by their degree fingerprints,
we may assume $F\equiv_2 H$.
Since $F\in\mathcal U_{n,p}$ and
$H\in\mathcal C_{n,p}$, the definition of
$\mathcal U_{n,p}$ gives
\[
|E(F)\triangle E(H)|
<\frac{2|E(F)|}{(\log n)^{1/4}}
\le\frac{2n^2p}{(\log n)^{1/4}}
=o(n^2p).
\]
This contradicts Theorem~\ref{thm:correction}
for all sufficiently large $n$.
Thus no two nonisomorphic graphs in
$\mathcal U_{n,p}$ have the same deck.
By Lemma~\ref{lem:fibre},
\[
\Pbb(G(n,p)\text{ is not reconstructible})
\le 2\Pbb(G(n,p)\notin\mathcal U_{n,p})
=o(1).
\]

Finally, \cite[Theorem~3]{JohnstonKronenbergRobertsScott} establishes
reconstructibility with high probability for
$p\ge n^{-2/3-\gamma}$ when $\gamma>0$
is sufficiently small.

\end{proof}

\appendix
\section{Proof of the preliminaries}\label{app:preliminaries}

\begin{proof}[Proof of Lemma~\ref{lem:basic}]
Fix an even integer $C>2(M+2)/(3\gamma)$.
Since $np\gg\log n$, Chernoff's inequality gives
\begin{equation}\label{eq:basic-degree}
\Pbb\bigl(\exists v:\deg(v)\notin[np/2,2np]\bigr)
\le 2ne^{-cnp}=o(n^{-M}).
\end{equation}

For distinct vertices $x,y$, their codegree has distribution
$\operatorname{Bin}(n-2,p^2)$, and hence
\begin{equation}\label{eq:basic-codegree}
\Pbb\bigl(|N(x)\cap N(y)|\ge C/2\bigr)
\le \left(\frac{2enp^2}{C}\right)^{C/2}.
\end{equation}

Conditional on $N(x)$ and $N(y)$, there are at most
$\deg(x)\deg(y)$ possible middle edges, independently
present with probability $p$, each contributing to at
most two 3-paths. Thus,
\begin{equation}\label{eq:basic-path}
\Pbb\bigl(P_3(x,y)\ge C,\
\deg(x),\deg(y)\le 2np\bigr)
\le \left(\frac{8en^2p^3}{C}\right)^{C/2}.
\end{equation}

Similarly, conditional on $N(x)$, the number $T_x$
of triangles containing $x$ satisfies
$T_x\sim\operatorname{Bin}(\binom{\deg(x)}{2},p)$.
Thus,
\begin{equation}\label{eq:basic-triangle}
\Pbb\bigl(T_x\ge C/2,\ \deg(x)\le 2np\bigr)
\le \left(\frac{4en^2p^3}{C}\right)^{C/2}.
\end{equation}

Since $np^2\le n^{-1/3-2\gamma}$ and
$n^2p^3\le n^{-3\gamma}$, combining
\eqref{eq:basic-degree}--\eqref{eq:basic-triangle}
with a union bound shows that, with probability
$1-O(n^{-M})$, all degrees lie in $[np/2,2np]$,
while all codegrees, 3-path counts, and triangle
counts are at most $C$.
\end{proof}

\begin{proof}[Proof of Corollary~\ref{lem:disjoint-fourpaths}]
Let $E_K$ be the event in
Lemma~\ref{lem:basic} hold with constant $K$.
By Lemma~\ref{lem:basic}, $\Pbb(E_K^c)=o(1)$. On $E_K$, for distinct vertices $x,y$ and
$z\notin\{x,y\}$, the number of 4-paths from $x$ to $y$
containing $z$ internally is at most
\begin{equation}\label{eq:fourpath-incidence}
\begin{aligned}
&P_3(z,y)+|N(x)\cap N(z)|\,|N(z)\cap N(y)|
+P_3(x,z)\\
&\qquad\le K^2+2K.
\end{aligned}
\end{equation}

Since each 4-path has three internal vertices, a greedy
selection from any family of $\kappa$ distinct 4-paths yields
at least $\left\lfloor\frac{\kappa}{3(K^2+2K)}\right\rfloor$ paths with pairwise disjoint internal vertex sets.
Taking $C=3(K^2+2K)$ completes the proof.
\end{proof}

\begin{proof}[Proof of Lemma~\ref{lem:fourpaths}]
Recall $E_K$, and set $C_0=3(K^2+2K)$. For distinct vertices $x,y$, let $F_{x,y,r}$ denote
the event that there exist $r$ 4-paths between $x$
and $y$ with pairwise disjoint internal vertex sets.
There are at most $n^{3r}/r!$ candidate families
of such paths, each requiring $4r$ distinct edges.
Hence,
\begin{equation}\label{eq:disjoint-fourpaths}
\Pbb(F_{x,y,r})
\le \frac{n^{3r}p^{4r}}{r!}
\le \left(\frac{en^3p^4}{r}\right)^r.
\end{equation}

Take $r=\lfloor 4en^3p^4\rfloor$ and $C=4eC_0$.
By Corollary~\ref{lem:disjoint-fourpaths}, on $E_K$, $\{P_4(x,y)>Cn^3p^4\}\subseteq F_{x,y,r}.$ Since $n^3p^4\ge(\log n)^{12}$, we have
$r\ge 2en^3p^4$ for all sufficiently large $n$.
Thus, by \eqref{eq:disjoint-fourpaths} and a union bound,
\begin{align*}
\Pbb\bigl(\exists x\ne y:
P_4(x,y)>Cn^3p^4\bigr)
&\le \Pbb(E_K^c)
+\sum_{x\ne y}\Pbb(F_{x,y,r})\\
&\le O(n^{-M})+n^2 2^{-r}
=O(n^{-M}),
\end{align*}
\end{proof}

\begin{proof}[Proof of Lemma~\ref{lem:cycles}]
Every 5-cycle containing an edge $uv$ corresponds to a
unique simple 4-path from $u$ to $v$. Hence,
\begin{equation}\label{eq:edge-cycle-sum}
c_5(G;uv)=P_4^G(u,v)
=\sum_{w\in N(u)\setminus\{v\}}P_3^{G-u}(v,w).
\end{equation}
The upper bound follows from Lemma~\ref{lem:fourpaths}.

For the lower bound, recall $E_K$. Fix distinct vertices
$u,v$, and define the sets of exclusive neighbors
of $u$ and $v$ by $A=N(u)\setminus(N(v)\cup\{v\}),$ $B=N(v)\setminus(N(u)\cup\{u\}).$ On $E_K$, the sets $A,B$ are disjoint and satisfy
$|A|,|B|\asymp np$.

Let $Z_{u,v}$ count vertices outside
$N(u)\cup N(v)\cup\{u,v\}$ having at least one
neighbor in each of $A$ and $B$. 
Conditional on $N(u)$ and $N(v)$,
\[
Z_{u,v}\sim\operatorname{Bin}\left(
n-|N(u)\cup N(v)\cup\{u,v\}|,
\bigl(1-(1-p)^{|A|}\bigr)
\bigl(1-(1-p)^{|B|}\bigr)
\right).
\]

Since $|A|,|B|\asymp np$ and $np^2=o(1)$,
\[
\begin{aligned}
\Ebb[Z_{u,v}\mid N(u),N(v)]
&=\bigl(n-|N(u)\cup N(v)\cup\{u,v\}|\bigr)\\
&\quad\times\bigl(1-(1-p)^{|A|}\bigr)
\bigl(1-(1-p)^{|B|}\bigr)\\
&\asymp n\cdot(|A|p)(|B|p)
\asymp n^3p^4.
\end{aligned}
\]

Thus, Chernoff's inequality gives constants $c,c'>0$
such that
\[
\Pbb\bigl(Z_{u,v}<cn^3p^4\mid N(u),N(v)\bigr)
\le e^{-c'n^3p^4}
\]
uniformly over these neighborhoods.

Each vertex counted by $Z_{u,v}$ yields a distinct
4-path from $u$ to $v$, so $P_4^G(u,v)\ge Z_{u,v}$.
By \eqref{eq:edge-cycle-sum} and a union bound,
\begin{align*}
&\Pbb\bigl(\exists uv\in E(G):
c_5(G;uv)<cn^3p^4\bigr)\le\Pbb(E_K^c)+n^2e^{-c'n^3p^4}
=o(1),
\end{align*}
where we used $n^3p^4\ge(\log n)^{12}$.

For two adjacent edges $uv,uw$, deleting $u$ from
any 5-cycle containing both leaves a 3-path from
$v$ to $w$. Hence, at most $P_3^{G-u}(v,w)\le C$
such cycles exist.

For two disjoint edges $ab,cd$, any 5-cycle
containing both must contain one of the edges
$ac,ad,bc,bd$. Once this single edge is fixed, the
remaining endpoints must be joined by a 2-path,
of which there are at most $C$ by the codegree limitation.
Thus, any two distinct edges lie together in at most
$4C$ 5-cycles, completing the proof.

\end{proof}

\begin{proof}[Proof of Lemma~\ref{lem:overlap}]
Recall $E_K$, and let $E_K'$ be its
intersection with the event in Lemma~\ref{lem:fourpaths}.
Then $\Pbb((E_K')^c)=o(1)$. 
Fix distinct vertices $u,v,z$, and let $G'=G-u$.
Write $B_{G'}(x,r)$ for the radius-$r$ ball
around $x$ in $G'$.
On $E_K'$, \eqref{eq:P4_bound} imply $P_5^{G'}(v,z)\le Cn^4p^5,$ $P_6^{G'}(v,z)\le Cn^5p^6,$ and $|B_{G'}(x,r)|=O((np)^r)$ for $0\le r\le2$.

For each $w\in B_{G'}(v,3)\cap B_{G'}(z,3)$,
consider a breadth-first-search tree rooted at $w$.
Let $a$ be the last common vertex of the tree paths
from $w$ to $v$ and $z$, and let $\ell$ be the length
of the resulting simple $v$--$z$ path.
Then $1\le\ell\le6$, and
\[
2\operatorname{dist}_{G'}(w,a)+\ell
=\operatorname{dist}_{G'}(w,v)
+\operatorname{dist}_{G'}(w,z)\le6.
\]
Consequently, $\operatorname{dist}_{G'}(w,a)
\le r_\ell:=3-\left\lceil\frac{\ell}{2}\right\rceil.$ Since $a$ lies on a simple $v$--$z$ path of length
$\ell$, the balls of radius $r_\ell$ centered at
the vertices of all such paths cover
$B_{G'}(v,3)\cap B_{G'}(z,3)$.

\begin{table}[htbp]
\centering
\caption{Number of paths and their contributions.}
\label{tab:overlap-paths}
\begin{tabular}{cccc}
\toprule
$\ell$ & $P_\ell^{G'}(v,z)$
& $r_\ell$ & Contribution \\
\midrule
$1$ & $\le 1$             & $2$ & $O(n^2p^2)$ \\
$2$ & $\le K$             & $2$ & $O(n^2p^2)$ \\
$3$ & $\le K$             & $1$ & $O(np)$ \\
$4$ & $\le Cn^3p^4$       & $1$ & $O(n^4p^5)$ \\
$5$ & $\le 2Cn^4p^5$      & $0$ & $O(n^4p^5)$ \\
$6$ & $\le 4Cn^5p^6$      & $0$ & $O(n^5p^6)$ \\
\bottomrule
\end{tabular}
\end{table}

Table~\ref{tab:overlap-paths} summarizes the bounds
on the number of paths and their contributions
to the covering estimate. Summing over the path lengths, we obtain
\begin{equation}\label{eq:overlap-balls}
\begin{aligned}
|B_{G'}(v,3)\cap B_{G'}(z,3)|
&\le \sum_{\ell=1}^6(\ell+1)P_\ell^{G'}(v,z)
\max_x|B_{G'}(x,r_\ell)|\\
&\le Cn^5p^6.
\end{aligned}
\end{equation}

By the definition of $\Gamma_u$,
\[
\Gamma_u(v)\cap\Gamma_u(z)
\subseteq N(u)\cap B_{G'}(v,3)\cap B_{G'}(z,3).
\]
Conditional on $G'$,
\[
|N(u)\cap B_{G'}(v,3)\cap B_{G'}(z,3)|
\sim\operatorname{Bin}\left(
|B_{G'}(v,3)\cap B_{G'}(z,3)|,p
\right).
\]
Whenever \eqref{eq:overlap-balls} holds,
Chernoff's inequality gives, for sufficiently
large $C$,
\[
\begin{aligned}
&\Pbb\bigl(
|N(u)\cap B_{G'}(v,3)\cap B_{G'}(z,3)|
>C(n^5p^7+\log n)\mid G'\bigr)\le n^{-5}.
\end{aligned}
\]

Since \eqref{eq:overlap-balls} depends only on $G'$
and holds on $E_K'$, averaging over $G'$ and
taking a union bound over all distinct $u,v,z$
yields
\[
\begin{aligned}
&\Pbb\bigl(\exists\,u,v,z\text{ distinct}:
|\Gamma_u(v)\cap\Gamma_u(z)|
>C(n^5p^7+\log n)\bigr)\le\Pbb((E_K')^c)+n^3n^{-5}=o(1).
\end{aligned}
\]

Finally,
\[
\frac{n^5p^7+\log n}{n^3p^4}
=n^2p^3+\frac{\log n}{n^3p^4}
\le n^{-3\gamma}+(\log n)^{-11}
=o(1),
\]
completing the proof.
\end{proof}

\begin{proof}[Proof of Lemma~\ref{lem:spectral}]
By the $(S,t)$-rope estimate of Baron and Kahn~\cite{BaronKahn},
with $t=4$ and $\beta=2s/(n^2p)$, the number of simple
4-paths whose terminal edges belong to $S$ is, with high
probability, uniformly over all $S\subseteq E(G)$, at most
\[
O\!\left(\max\{\beta^2n^5p^4,\,
\beta n^4p^3\}\right)
=O(np^2s^2+n^2p^2s).
\]
Each unoriented 4-path corresponds to exactly two ordered
vertex sequences in \eqref{eq:rope-count}, and the factor
of two is absorbed into $C$.
\end{proof}

\begin{proof}[Proof of Lemma~\ref{lem:links}]
Recall $E_K$. 
Fix a vertex $u$, and let $\mathcal F_u$ be the
$\sigma$-field generated by all edges incident to
$u$ or a vertex in $N(u)$.
For $J\subseteq N(u)$ with $j=|J|\ge1$, let
$Z_J$ count edges between vertices outside
$N(u)\cup\{u\}$ that each have a neighbor in $J$.
On $E_K$, at most
$\sum_{v\in J}\deg(v)\le 2npj$ vertices
have a neighbor in $J$. Thus, conditional on
$\mathcal F_u$,
\[
Z_J\preceq\operatorname{Bin}
\left(\lceil 2n^2p^2j^2\rceil,p\right),
\]
where $\preceq$ denotes stochastic domination.

By Lemma~\ref{lem:basic}, each endpoint of an edge
counted by $Z_J$ has at most $K$ neighbors in $J$.
Thus, each such edge witnesses at most $K^2$
edges of $Q_u(G)[J]$. The remaining witnessing 3-paths have at least
one internal vertex in $N(u)$. Orient such a path
as $v-x-y-w$ with $x\in N(u)$.
There are at most $2e(G[N(u)])\le2K$ choices
for $(v,x)$ and at most $P_3^G(u,x)\le K$
choices for the remaining path. Hence,
\begin{equation}\label{eq:small-link-decomposition}
e(Q_u(G)[J])\le K^2Z_J+2K^2.
\end{equation}

\begin{figure}[H]
    \centering
    \includegraphics[width=0.60\linewidth]
    {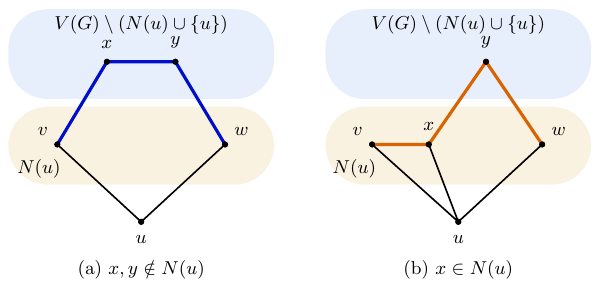}
    
\caption{Two configurations of a 3-path $v-x-y-w$
witnessing an edge of $Q_u(G)[J]$, where
$v,w\in J\subseteq N(u)$.
(a) Both internal vertices lie outside $N(u)$,
and the middle edge $xy$ is counted by $Z_J$.
(b) The vertex $x$ lies in $N(u)$, so
$vx\in E(G[N(u)])$ and $u-w-y-x$ is a 3-path.}
\label{fig:small-link-witnesses}

\label{fig:correction-cases11}
    \label{fig:correction-cases}
\end{figure}

Choose
$\delta_\eta=\min\{1,\eta/(16eK^2)\}$.
For $j\le\delta_\eta np$, the binomial tail bound
gives
\begin{align*}
&\Pbb\left(
Z_J>\frac{\eta n^3p^4j}{2K^2}
\,\middle|\,\mathcal F_u\right)\le
\left(\frac{4eK^2\delta_\eta}{\eta}\right)^{
\eta n^3p^4j/(2K^2)}
\le e^{-c_{\eta,K}n^3p^4j}.
\end{align*}

Since $n^3p^4\to\infty$, we also have
$2K^2\le\eta n^3p^4j/2$ for sufficiently
large $n$. By \eqref{eq:small-link-decomposition},
\begin{equation}\label{eq:small-link-tail}
\Pbb\bigl(
e(Q_u(G)[J])>\eta n^3p^4j
\mid\mathcal F_u\bigr)
\le e^{-c_{\eta,K}n^3p^4j}.
\end{equation}

Finally, since $|N(u)|\le2np$ on $E_K$
and $n^3p^4\gg\log n$, a union bound over
all roots and eligible sets gives
\begin{align*}
&\Pbb\bigl(\exists u,\ J\subseteq N(u):
1\le |J|\le\delta_\eta np,\
e(Q_u(G)[J])>\eta n^3p^4|J|\bigr)\\
&\quad\le\Pbb(E_K^c)
+n\sum_{j=1}^{\lfloor\delta_\eta np\rfloor}
\binom{\lfloor2np\rfloor}{j}
e^{-c_{\eta,K}n^3p^4j}
=o(1).
\end{align*}
\end{proof}

\begin{proof}[Proof of Proposition~\ref{prop:single}]
Let $\mathcal E$ be the event on which the conclusions
of Lemmas~\ref{lem:basic}, \ref{lem:fourpaths}
and~\ref{lem:cycles}, and Corollary~\ref{lem:disjoint-fourpaths}
hold. Then $\Pbb(\mathcal E)=1-o(1)$.

Suppose that $G\in\mathcal E$ admits a counterexample
$B\subseteq E(G)$ with
$0<|B|\le\rho_\varepsilon n^2p$ and
$J_G(B)>\varepsilon n^3p^4|B|$.
For each $e\in B$, let $r_e$ denote the number of
5-cycles counted at $e$ in $J_G(B)$.
Since $r_e\le Cn^3p^4$ by Lemma~\ref{lem:cycles},
at least $c_\varepsilon|B|$ edges satisfy $r_e\ge\frac{\varepsilon}{2}n^3p^4,$ where $c_\varepsilon=\varepsilon/(2C)>0$.
Call these edges active.

Deleting an active edge $e$ from each corresponding
5-cycle leaves a 4-path between its endpoints containing
an internal $B$-edge. 
By Corollary~\ref{lem:disjoint-fourpaths}, for each
active edge $e$, we can choose a family $\mathcal P_e$
of at least $c'_\varepsilon n^3p^4$ such paths with
pairwise disjoint internal vertex sets, for some
fixed $c'_\varepsilon>0$.
Since each path contains an internal $B$-edge,
these edges are distinct across the family, giving
$|B|\ge c'_\varepsilon n^3p^4$.
For each graph admitting a counterexample,
fix one such $B$ and the families $\mathcal P_e$
deterministically.

Choose $0<\theta\le\min\{1/4,c_\varepsilon/480\}$
and retain each edge of $G$ independently with
probability $\theta$, obtaining $G_0\sim G(n,p_0)$,
where $p_0=\theta p$.
Let $S=B\cap E(G_0)$ and $s=|S|$. Since the paths in each $\mathcal P_e$ have disjoint
edge sets, they survive independently with probability $\theta^4$. For every active $e$,
by Chernoff's inequality and a union bound,
uniformly over the chosen counterexamples,
\begin{equation}
\begin{aligned}
&\Pbb\Bigl(
\frac{\theta|B|}{2}\le s\le2\theta|B|,\ 
\#\{P\in\mathcal P_e:P\subseteq G_0\}
\ge\frac{c'_\varepsilon}{2}n^3p_0^4
\,\Bigm|\,G\Bigr)\\
&\quad\ge 1-
n^2e^{-c'_\varepsilon n^3p_0^4/8}
-2e^{-\theta|B|/8}
=1-o(1),
\end{aligned}
\label{eq:conditional_prob_lower_bound}
\end{equation}
where we used $|B|\ge c'_\varepsilon n^3p^4$
and $n^3p^4\gg\log n$.


For $S\subseteq E(G_0)$, let $R(S)$ be the set of
unordered vertex pairs joined by at least
$(c'_\varepsilon/2)n^3p_0^4$ 4-paths in $G_0$
with pairwise disjoint internal vertex sets,
each containing an internal $S$-edge. Whenever \eqref{eq:conditional_prob_lower_bound} holds, the endpoints of every
active edge absent from $G_0$ belong to $R(S)$.
Since $c_\varepsilon|B|-s\ge c_\varepsilon|B|/2$
and $s\le2\theta|B|$, we obtain
\begin{equation}\label{eq:upward-event}
|(E(G)\setminus E(G_0))\cap R(S)|
\ge\frac{c_\varepsilon s}{4\theta},
\qquad
1\le s\le2\rho_\varepsilon n^2p_0.
\end{equation}

We next show that \eqref{eq:upward-event} is unlikely.
Let $\mathcal E_0$ be the event given by
Lemmas~\ref{lem:basic}, \ref{lem:fourpaths}
and~\ref{lem:spectral} at density $p_0$.
Since $\theta$ is fixed, the same bounds apply
at $p_0$, and $\Pbb(\mathcal E_0)=1-o(1)$. 

Fix $G_0\in\mathcal E_0$ and
$S\subseteq E(G_0)$ with $s\ge1$.
For each pair $\{x,y\}\in R(S)$, choose a qualifying
family of 4-paths and one internal $S$-edge from
each path. Orient the paths consistently so that
at least half have their chosen $S$-edge in the
second position.

Any two such paths $x-a_1-b_1-c_1-y,$ and $x-a_2-b_2-c_2-y$ give rise to two simple 4-paths $b_1-a_1-x-a_2-b_2,$ $b_1-c_1-y-c_2-b_2.$ The first has both terminal edges in $S$ and
is therefore an $S$-rope.
Conversely, these two paths determine the original
pair $\{x,y\}$ and the chosen paths, up to a
constant number of choices.

\begin{figure}[H]
    \centering
    \includegraphics[width=0.5\linewidth]
    {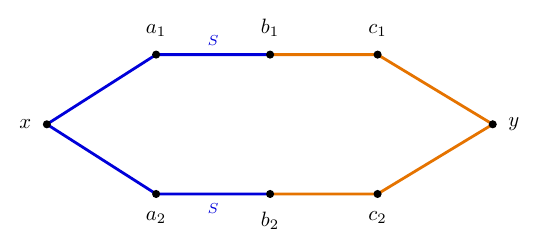}

\caption{Two internally 4-paths
from $x$ to $y$, with $a_1b_1,a_2b_2\in S$.
The blue path $b_1-a_1-x-a_2-b_2$ is an $S$-rope,
whose two terminal edges belong to $S$.}

\label{fig:small-link-witnesses}

\label{fig:correction-cases11}
    \label{fig:correction-cases}
\end{figure}

Each pair in $R(S)$ thus gives at least $cn^6p_0^8$
pairs of paths, while each $S$-rope admits at most
$Cn^3p_0^4$ choices for the second 4-path by
Lemma~\ref{lem:fourpaths}. Applying
Lemma~\ref{lem:spectral} to count the $S$-ropes,
we obtain
\begin{equation}\label{eq:candidate-bound}
|R(S)|\le C\frac{np_0^2s^2+n^2p_0^2s}{n^3p_0^4}
=\frac{Cs}{p_0}\left(\frac{s}{n^2p_0}+\frac{1}{np_0}\right).
\end{equation}

Since $s\ge1$ and $p_0=\theta p$, the lower bound
on $p$ implies
$\frac{1}{np_0}\le(n^2p_0)^{-1/6}
\le(s/(n^2p_0))^{1/6}$
for all sufficiently large $n$.
Thus, for $s\le n^2p_0$,
\eqref{eq:candidate-bound} yields
\[
|R(S)|\le\frac{Cs}{p_0}
\left(\frac{s}{n^2p_0}\right)^{1/6}.
\]

Conditional on $G_0$, each nonedge of $G_0$ is
added independently to form $G$ with probability
$\frac{p-p_0}{1-p_0}\le p$.
Hence, for
$X_S=|(E(G)\setminus E(G_0))\cap R(S)|$,
\[
\Ebb[X_S\mid G_0]
\le\frac{Cs}{\theta}
\left(\frac{s}{n^2p_0}\right)^{1/6}.
\]

Choose $\rho_\varepsilon>0$ sufficiently small that
$2\rho_\varepsilon<e^{-1}$ and $\frac{4eC}{c_\varepsilon}(2\rho_\varepsilon)^{1/12}\le1.$ For $1\le s\le2\rho_\varepsilon n^2p_0$,
the binomial tail bound gives
\begin{equation}\label{eq:upward-tail}
\begin{aligned}
\Pbb\left(
X_S\ge\frac{c_\varepsilon s}{4\theta}
\,\middle|\,G_0\right)
&\le\left(
\frac{4eC}{c_\varepsilon}
\left(\frac{s}{n^2p_0}\right)^{1/6}
\right)^{c_\varepsilon s/(4\theta)}\\
&\le\left(\frac{s}{n^2p_0}\right)^{
c_\varepsilon s/(48\theta)}
\le\left(\frac{s}{n^2p_0}\right)^{10s}.
\end{aligned}
\end{equation}

Let $\mathcal E_{\mathrm{up}}$ be the event that
\eqref{eq:upward-event} holds for some
$S\subseteq E(G_0)$.
For every $G_0\in\mathcal E_0$, we have
$|E(G_0)|\le n^2p_0$. Thus, by
\eqref{eq:upward-tail} and a union bound,
\begin{equation}\label{eq:upward-union}
\begin{aligned}
\Pbb(\mathcal E_{\mathrm{up}}\mid G_0)
&\le\sum_{1\le s\le2\rho_\varepsilon n^2p_0}
\binom{|E(G_0)|}{s}
\left(\frac{s}{n^2p_0}\right)^{10s}\\
&\le\sum_{1\le s\le2\rho_\varepsilon n^2p_0}
\left(\frac{s}{n^2p_0}\right)^{8s}
=o(1).
\end{aligned}
\end{equation}

By \eqref{eq:conditional_prob_lower_bound},  $\Pbb(\mathcal E_{\mathrm{up}}\mid G)\ge1-o(1)$
uniformly over all $G\in\mathcal E$ admitting
a counterexample. On the other hand,
\eqref{eq:upward-union} and
$\Pbb(\mathcal E_0^c)=o(1)$ imply
$\Pbb(\mathcal E_{\mathrm{up}})=o(1)$.
Consequently, the probability that
$G\in\mathcal E$ admits a counterexample is $o(1)$.
Since $\Pbb(\mathcal E^c)=o(1)$,
the proof is complete.

\end{proof}

\section{Proof of
Proposition~\ref{prop:moved}}\label{app:correction}

This section is devoted to the proof of
Proposition~\ref{prop:moved}. The key is to condition
on the neighborhoods of the source and target vertices,
leaving independent edges between their private neighbors.
We begin with a uniform bound on the clique number
of the auxiliary graph defined by 3-paths.

\begin{lemma}\label{lem:clique}
There is a fixed integer $R_*$ such that, with high probability,
\[
 \omega\bigl(\mathcal P_3(G)\bigr)\le R_*,
\]
where $\omega$ denotes the maximum size of a clique, and $\mathcal P_3(F)$ joins distinct vertices whenever a simple 3-path joins them in $F$.
\end{lemma}

\begin{proof}
Fix an integer $d\ge1$, and choose $k_d$ sufficiently
large that $k_d>2d^2+1,$ and
\[
\frac{3(k_d-1)}{2(k_d-1+d^2)}
>\frac{1}{2/3+\gamma}.
\]
We first claim that, with high probability, no
independent set $I$ of size $k_d$ admits a choice
of one simple 3-path between each pair of its vertices
such that every internal vertex has at most $d$
neighbors in $I$.

Suppose such paths exist, and let $\widetilde G$
be their union. Write $L$ for the number of distinct
middle edges of these paths. Since $I$ is independent,
all internal vertices lie outside $I$.
Each middle edge belongs to at most $d^2$ chosen
paths, while every internal vertex is incident with
a middle edge and an edge joining it to $I$.
Hence, $L\ge\frac{k_d(k_d-1)}{2d^2}>k_d,$ $|V(\widetilde G)|-k_d\le2L,$
and
\[
|E(\widetilde G)|
\ge L+|V(\widetilde G)|-k_d.
\]
Consequently,
\begin{equation}\label{eq:path-union-density}
\begin{aligned}
\frac{|E(\widetilde G)|}{|V(\widetilde G)|}
&\ge\frac{L+|V(\widetilde G)|-k_d}
{|V(\widetilde G)|}\ge\frac{3L}{k_d+2L}\ge\frac{3(k_d-1)}{2(k_d-1+d^2)}
>\frac1{2/3+\gamma}.
\end{aligned}
\end{equation}

For each fixed $d\ge1$, define
\[
\mathcal E_d^{\mathrm{dens}}
=\left\{
\frac{|E(H)|}{|V(H)|}\le\frac1{2/3+\gamma}
\text{ for every }H\subseteq G
\text{ with }1\le|V(H)|\le k_d^2
\right\}.
\]
By \eqref{eq:path-union-density} and a union bound,
\[
\begin{aligned}
\Pbb\bigl((\mathcal E_d^{\mathrm{dens}})^c\bigr)
&\le\sum_{\substack{
H:\,1\le|V(H)|\le k_d^2\\
|E(H)|/|V(H)|>1/(2/3+\gamma)}}
n^{|V(H)|}p^{|E(H)|}
=o(1),
\end{aligned}
\]
the claim holds on
$\mathcal E_d^{\mathrm{dens}}$.

Recall $E_K$, with $K$ chosen as a positive integer,
and work on $E_K\cap\mathcal E_K^{\mathrm{dens}}$ with $d=K$.
Fix a clique $W$ of $\mathcal P_3(G)$ and a vertex $a$.
Let $W_a$ be obtained from $W\cap N_G(a)$ by
removing all vertices belonging to a triangle
containing $a$. At most $2K$ vertices are removed,
and $W_a$ is independent.

Moreover, a simple 3-path between two vertices
of $W_a$ cannot contain $a$ internally, since
otherwise one endpoint would lie in a triangle
containing $a$. Thus, every internal vertex $x$
of such a path satisfies
\[
|N_G(x)\cap W_a|
\le |N_G(x)\cap N_G(a)|\le K.
\]

On $\mathcal E_K^{\mathrm{dens}}$, we have
$|W_a|<k_K$.
 Therefore, for every $a\in V(G)$,
\begin{equation}\label{eq:clique-root-neighbours}
|W\cap N_G(a)|\le 2K+k_K-1.
\end{equation}

By \eqref{eq:clique-root-neighbours}, $G[W]$ has
maximum degree at most $2K+k_K-1$, and hence
contains an independent set $I'\subseteq W$ with
$|I'|\ge |W|/(2K+k_K)$.
Since every vertex of $G$ has at most
$2K+k_K-1$ neighbors in $I'$, and $I'$ remains
a clique of $\mathcal P_3(G)$, the event
$\mathcal E_{2K+k_K-1}^{\mathrm{dens}}$ implies
$|I'|<k_{2K+k_K-1}$. Consequently,
\[
|W|\le (2K+k_K)(k_{2K+k_K-1}-1)=:R_*.
\]
Since $E_K$, $\mathcal E_K^{\mathrm{dens}}$ and
$\mathcal E_{2K+k_K-1}^{\mathrm{dens}}$ all hold
with high probability, the proof is complete.

\end{proof}

\begin{proof}[Proof of Proposition~\ref{prop:moved}]

Let $c_0>0$ be the constant
obtained from Lemmas~\ref{lem:basic} and~\ref{lem:cycles},
and choose $\delta>0$ from Lemma~\ref{lem:links}
with $\eta=c_0/8$.
Recall $R_*$ is the constant in Lemma~\ref{lem:clique},
and fix $0<\eps\le\min\{1/4,\delta/(192R_*)\}$. Let $\mathcal E_{\mathrm{mov}}$ be the intersection
of the events in Lemmas~\ref{lem:basic},
\ref{lem:cycles}, \ref{lem:links},
\ref{lem:overlap} and~\ref{lem:clique}. For each fixed nonempty $D$, the number of vertices
outside $D$ having exactly one neighbor in $D$
has distribution $\operatorname{Bin}\left(
n-|D|,\ |D|p(1-p)^{|D|-1}
\right),$ with mean $(1-o(1))|D|np$, uniformly over
$|D|\le4np$, since $np^2=o(1)$. Chernoff's inequality and a union bound over
at most $n^{|D|}$ sets of each size show that $\Pbb(\mathcal E_{\mathrm{mov}})=1-o(1).$

Fix $u,T$ and a bijection $\phi:T\to N_G(u)$.
Let
\[
M=\{v\in N_G(u):\phi^{-1}(v)\ne v\}.
\]

If $M=\varnothing$, \eqref{eq:moved-loss} is immediate.
Suppose that $0<|M|\le\delta np$. 
By Lemmas~\ref{lem:cycles} and~\ref{lem:links},
the number of edges of $Q_u(G)$ between
$M$ and $N_G(u)\setminus M$ is at least
\[
c_0n^3p^4|M|-2e(Q_u(G)[M])
\ge\frac{3c_0}{4}n^3p^4|M|.
\]

For $v\in M$ and $w\in N_G(u)\setminus M$,
we have $\phi^{-1}(w)=w$.
If $vw\in E(Q_u(G))$ and
$\phi^{-1}(v)w\in E(\mathcal L_u^G(T))$,
then $w\in\Gamma_u^G(v)\cap
\Gamma_u^G(\phi^{-1}(v)).$ By Lemma~\ref{lem:overlap}, for each $v\in M$,
there are at most $o(n^3p^4)$ such vertices $w$,
uniformly over $u,T,\phi$.
Consequently, at least
$(c_0/2)n^3p^4|M|$ link edges are lost,
proving \eqref{eq:moved-loss} in this case.

Now suppose that $|M|>\delta np$.
For each $x\in N_G(u)\cup T$, let $A(x)$ be the set
of vertices outside $N_G(u)\cup T\cup\{u\}$
whose unique neighbor in $N_G(u)\cup T$ is $x$.
Since $|N_G(u)\cup T|\le4np$, on
$\mathcal E_{\mathrm{mov}}$ we have $|A(x)|\le(1+\eps)np$ and
\begin{equation}\label{eq:private-set-bounds}
\sum_{x\in N_G(u)\cup T}|A(x)|
\ge(1-2\eps)|N_G(u)\cup T|np.
\end{equation}
Consequently, at most
$6\eps|N_G(u)\cup T|\le24\eps np$
vertices $x\in N_G(u)\cup T$ satisfy
$|A(x)|<np/2$.
Call the remaining vertices good.

Fix $u$ and condition on $N_G(u)$. Fix any $T'\subseteq[n]\setminus\{u\}$ with
$|T'|=|N_G(u)|$ and a bijection
$\phi':T'\to N_G(u)$.
Condition further on all edges incident to
$N_G(u)\cup T'$, and then on all remaining edges
except those joining distinct private sets.
Denote the resulting information by $\mathcal F$.

For each unordered pair $P=\{x,y\}\subseteq
N_G(u)\cup T'$, let $Z_P$ indicate the presence
of an edge between $A(x)$ and $A(y)$.
Conditional on $\mathcal F$, the variables
$Z_P$ are independent, with
\[
\Pbb(Z_P=1\mid\mathcal F)
=1-(1-p)^{|A(x)||A(y)|}.
\]

By \eqref{eq:private-set-bounds} and the definition
of good vertices, since $n^2p^3=o(1)$, we have
\begin{equation}\label{eq:private-pair-probabilities}
\begin{aligned}
\Pbb(Z_P=1\mid\mathcal F)&\le Cn^2p^3,\\
\Pbb(Z_P=1\mid\mathcal F)&\ge cn^2p^3.
\end{aligned}
\end{equation}
The second bound holds whenever both vertices
of $P$ are good.

Let $Q_P$ indicate the existence
of a 3-path between $x$ and $y$ in $G-u$,
and let $Q_P^{(0)}$ be such a path using
only edges determined by $\mathcal F$.
Since any unconditioned edge on such a path
must be its middle edge joining $A(x)$ and $A(y)$,
we have
\begin{equation}\label{eq:background-private-decomposition}
Q_P=Q_P^{(0)}\vee Z_P.
\end{equation}

By Lemma~\ref{lem:clique}, the graph on
$N_G(u)\cup T'$ with edges $\{P:Q_P^{(0)}=1\}$
has clique number at most $R_*$.
Moreover, its restriction to $N_G(u)$ has
maximum degree at most $Cn^3p^4$ by
Lemma~\ref{lem:cycles}. Let $M'=\{v\in N_G(u):(\phi')^{-1}(v)\ne v\},$ and suppose $|M'|>\delta np$. 
By Tur\'an's theorem~\cite{Turan}, there are at least
$|M'|^2/(2R_*)-|M'|/2$ pairs $P\subseteq M'$
satisfying $Q_{(\phi')^{-1}(P)}^{(0)}=0$.

Exclude pairs $P\subseteq M'$ for which
$(\phi')^{-1}(P)=P$, $Q_P^{(0)}=1$, or
$P$ contains a vertex that is not good.
At most $|M'|/2$ pairs satisfy the first condition,
since each corresponds to a transposition;
at most $C|M'|n^3p^4/2$ satisfy the second;
and at most $24\eps np|M'|$ satisfy the third.
Consequently, at least
\[
\frac{|M'|^2}{2R_*}-|M'|
-\frac{C|M'|n^3p^4}{2}-24\eps np|M'|
\ge\frac{|M'|^2}{4R_*}
\]
pairs remain, where we used $|M'|>\delta np$,
$n^2p^3=o(1)$, and
$\eps\le\delta/(192R_*)$.

For each remaining pair $P$, we have
$(\phi')^{-1}(P)\ne P$,
$Q_P^{(0)}=Q_{(\phi')^{-1}(P)}^{(0)}=0$,
and both vertices of $P$ are good.

Form a directed graph whose vertices are unordered
pairs, with an edge $P\to(\phi')^{-1}(P)$ for each
remaining $P$. Since $\phi'$ is a bijection,
every vertex has indegree and outdegree at most one.
A greedy selection therefore yields at least
$|M'|^2/(12R_*)$ edges with pairwise disjoint
endpoints.

By \eqref{eq:background-private-decomposition}
and \eqref{eq:private-pair-probabilities},
each selected comparison satisfies
\[
\begin{aligned}
&\Pbb\bigl(
Q_P=1,\ Q_{(\phi')^{-1}(P)}=0
\mid\mathcal F\bigr)=
\Pbb(Z_P=1\mid\mathcal F)
\Pbb(Z_{(\phi')^{-1}(P)}=0\mid\mathcal F)
\ge cn^2p^3.
\end{aligned}
\]

Since the selected comparisons are conditionally
independent and $|M'|^2n^2p^3\ge\delta^2n^4p^5$,
Chernoff's inequality gives
\[
\Pbb\left(
\sum_{\text{selected }P}
\mathbf1_{\{Q_P=1,\,
Q_{(\phi')^{-1}(P)}=0\}}
<c_\delta n^4p^5
\,\middle|\,\mathcal F
\right)
\le \exp(-c_\delta n^4p^5).
\]

Since there are at most $\exp(Cnp\log n)$ choices
of $u,N_G(u),T,\phi$, and $n^3p^4\gg\log n$,
a union bound gives
\begin{equation}\label{eq:large-moved-tail}
\begin{aligned}
&\Pbb\Bigl(\exists\,u,T,\phi:\ |M|>\delta np,\
\bigl|\phi^{-1}E(Q_u(G))
\setminus E(\mathcal L_u^G(T))\bigr|
<c_\delta n^4p^5\Bigr)\\
&\qquad\le\Pbb(\mathcal E_{\mathrm{mov}}^c)
+\exp(Cnp\log n-c_\delta n^4p^5)
=o(1).
\end{aligned}
\end{equation}
Since $|M|\le2np$, with high probability,
\[
\bigl|\phi^{-1}E(Q_u(G))
\setminus E(\mathcal L_u^G(T))\bigr|
\ge c_\delta n^4p^5
\ge\frac{c_\delta}{2}n^3p^4|M|
\]
uniformly over all $u,T,\phi$ with $|M|>\delta np$.

Combining the two cases with
$c=\min\{c_0/2,c_\delta/2\}$
proves \eqref{eq:moved-loss}.

\end{proof}

\section*{AI Declaration}

During the research and preparation of this manuscript,
the author used ChatGPT 6 Pro for mathematical and
editorial assistance. ChatGPT helped identify a classical
result of Baron and Kahn~\cite{BaronKahn}, which was used
to replace the original proof of Lemma~\ref{lem:spectral},
and assisted in simplifying the proof of
Proposition~\ref{prop:single}. Based on the author's strategy, ChatGPT also
helped develop and complete the proof of
Proposition~\ref{prop:moved}. Subsequently, the author
used ChatGPT to check mathematical arguments, prepare
figures, and polish
the manuscript. The author takes full responsibility for the
correctness of all results and proofs in this paper.

\section*{Acknowledgments}
The author is grateful to Shuyang Gong and Heng Ma
for generously sharing this problem and their progress
on it, and for many stimulating discussions that
played a crucial role in developing the rigidity
argument based on 5-cycle constraints.
 This work was supported by the National Natural Science
Foundation of China(No.~12595294; No.~12231002).

\end{document}